
\input amstex
\documentstyle{amsppt}
\loadmsbm
\UseAMSsymbols
\magnification=\magstep1
\vsize 7.2in
\font\rom cmr10

\rom
\baselineskip = 15 pt
\nopagenumbers
\def\lcm{\hbox{\rm lcm}}
\def\e{\varepsilon}

\topmatter
\title Residue classes free of values of Euler's function
\endtitle
\dedicatory Dedicated to Andrzej Schinzel on his sixtieth birthday
\enddedicatory
\author Kevin Ford, Sergei Konyagin and Carl Pomerance \endauthor
\address Department of Mathematics, Universtiy of Texas,
Austin, TX 78712, USA \endaddress
\address Department of Mechanics and Mathematics, Moscow State University,
Moscow 119899, Russia \endaddress
\address Department of Mathematics, The University of Georgia,
Athens, GA 30602, USA \endaddress
\endtopmatter

\document

{\bf 1. Introduction}
\bigskip
By a totient we mean a value taken by Euler's function $\phi(n)$.
Dence and Pomerance \cite{DP} have established
\proclaim{Theorem A} If a residue class contains
at least one multiple of 4, then it contains infinitely many totients.
\endproclaim
Since 1 is the only odd totient, it remains to examine residue classes
consisting entirely of numbers $\equiv 2\pmod{4}$.  In this paper we
shall characterize which of these residue classes contain infinitely
many totients and which do not.  We show that the union of all residue
classes that are totient-free has asymptotic density $3/4$, that is,
almost all numbers that are $\equiv2\pmod4$ are in a residue class
that is totient-free.  In the other direction, we show the existence
of a positive density of odd numbers $m$, such that for any $s\ge0$
and any even number $a$, the residue class $a\pmod{2^sm}$ contains 
infinitely many totients.

We remark that if a residue class $r\pmod{s}$ contains infinitely
many totients, it is possible, using the methods of \cite{DP} and
Narkiewicz \cite{N}, 
to get an asymptotic formula for the number of $n\le x$
with $\phi(n)\equiv r\pmod{s}$.
\smallskip
\noindent{\bf Acknowledgements.}  We take this opportunity to
thank Sybilla Beckmann-Kazez, Andrew Granville, Robert Rumely,
Andrzej Schinzel, Roy Smith, Robert Varley and Felipe Voloch
for helpful discussions.
Research of the second author was supported by 
Grant 96-01-00378 from the Russian
Foundation for Basic Research.
The third author is supported in part by an NSF grant.
\bigskip

{\bf 2.  Preliminary results}
\bigskip

Totients in a residue class consisting of numbers that are
$\equiv2\pmod4$ 
necessarily are of the form
$p^k-p^{k-1}$ for some prime $p\equiv 3\pmod{4}$ and $k\ge 1$.
We begin by characterizing those residue classes which contain only finitely
many totients.

\proclaim{Lemma 1}
Suppose $s\ge 1$, $k\ge 1$, $a\equiv 2\pmod 4$.  Then there is a number
$y\equiv 3\pmod{4}$ such that $y^k-y^{k-1} \equiv a \pmod{2^s}$.
\endproclaim

\demo{Proof}  The lemma is trivial when $s=1$ or $k=1$, so suppose $s\ge 2$,
$k\ge 2$.  It suffices to show that the congruence
$$
y^k-y^{k-1} \equiv x^k-x^{k-1} \pmod{2^s}
$$
has no solutions with $y,x\equiv 3\pmod{4}$ and $x\not\equiv y\pmod{2^s}$.
If such a solution exists, write $x=zy$, so that
$y(1-z^k) \equiv 1-z^{k-1} \pmod{2^s}$.  Since $z\not\equiv 1\pmod{2^s}$,
we have
$$
y(1+z+\cdots + z^{k-1}) \equiv 1+ \cdots + z^{k-2} \pmod{2^s}.
$$
However, as $y$ and $z$ are both odd,
the above congruence is impossible.
\qed\enddemo

\proclaim{Lemma 2} Suppose $k\ge 2$, $M\ge 1$ and $p\equiv 3\pmod{4}$ is prime.
Then there is a number $x$ with $(x,M)=1$ and
$x^k-x^{k-1}\equiv p^k-p^{k-1}\pmod{M}$.
\endproclaim

\demo{Proof}
It is sufficient to prove the
existence of such $x$ for $M=r^l$ where $r$ is a prime. If $r\neq p$
we set $x=p$.  If $r=p$ we look for $x=p^{k-1}u+1$ for some number $u$.
Then
$$
u(p^{k-1}u+1)^{k-1} \equiv p-1 \pmod{p^w}, \tag1
$$
where $w=\max(0,l-k+1)$.  Let
$$
f(U) = U(p^{k-1}U+1)^{k-1} -p + 1.
$$
Since $f(U) \equiv U+1 \pmod{p}$, which has the root $-1$,
and $f'(-1)\equiv1\pmod{p}$, Hensel's lemma implies there is some
root $u$ of (1).
\qed\enddemo

\proclaim{Lemma 3}  Suppose $m$ is odd, $s\ge 2$, $a\equiv 2\pmod4$.
If the congruence
$$
x^k-x^{k-1}\equiv a\pmod{m} \tag2
$$
has a solution with $k\ge 1$ and $(x,m)=1$, then the 
progression \hbox{$a\pmod{2^s m}$} contains infinitely many totients.
Otherwise the progression contains either one or no totients, according
as $a=p-1$ for some $p|m$ or not.

\endproclaim

\demo{Proof} 
Assume that (2) has such a solution.
By Lemma 1, there is a number $y\equiv 3\pmod4$ such
that $y^k-y^{k-1}\equiv a\pmod{2^s}$.  It follows from
Dirichlet's Theorem that there are infinitely many primes
$p\equiv x\pmod{m}$, $p\equiv y\pmod{2^s}$,
and for each we have $\phi(p^k) \equiv a\pmod{2^s m}$.

If (2) has no solution with $(x,m)=1$, the only possible solutions of
$\phi(z)\equiv a\pmod{2^sm}$ are $z=4$, $z=p^k$ or $z=2p^k$ where 
$p$ is an odd prime dividing $m$.  If $z=4$, then $a=2$,
implying (2) has the solution $x=2,k=2$, a contradiction.
In addition, by Lemma 2, if
$a\equiv p^k-p^{k-1} \pmod{m}$ for some odd prime $p$ and $k\ge 2$, then
(2) has a solution with $(x,m)=1$.  Hence $z$ is either a prime or twice
a prime dividing $m$.
\qed\enddemo

Using Lemma 3, it is possible to find residue classes consisting of
even numbers which are free of totients.  For example, the progressions
$302 \pmod{1092}$ and $790 \pmod{1092}$ contain no totients.  In verifying
this, since $1092=4\times3\times7\times13$, one only needs to check
(2) for $k$ up to 12.

In the other direction, we prove

\proclaim{Theorem 1}  Suppose $M=2^s m$, where $s\ge 2$ and $m$ is odd.
If $a=\phi(b)>1$, where $b$ is neither prime nor
twice an odd prime, then any arithmetic progression $a\pmod{M}$
contains infinitely many totients.
\endproclaim

\demo{Proof} If $a$ is divisible by 4, the result follows from
Theorem A. Otherwise $a=2$ or $a=p^k-p^{k-1}$ where $p$ is an odd 
prime, $k>1$.

If $a=2$, $M=2^sm$, $m$ is odd, then for any prime $q$ such that
$q\equiv -1\pmod{2^s}$, $q\equiv 2\pmod{m}$ we have
$\phi(q^2)\equiv 2\pmod{M}$.

In the case $a=p^k-p^{k-1}$,
by Lemma 2 there is an $x$ such that $(x,M)=1$
and $x^k-x^{k-1}\equiv a\pmod{M}$. For any prime $q\equiv x\pmod{M}$
we have $\phi(q^k)\equiv a\pmod{M}$. 
\qed\enddemo

\proclaim{Question}
Suppose $a\equiv 2\pmod{4}$
is either a non-totient or a totient with exactly two pre-images
$\{ p, 2p \}$ for some prime $p$.  Is $a$ contained in
a residue class containing no totients other than $a$ itself?
\endproclaim

The numbers 10 and 14 are the two smallest such $a$.
A short search using a computer reveals
that the progression
\hbox{$14 \pmod{2^2\times3\times5\times13\times37}$} contains no totients
and the progression
$$
10 \pmod{4M}, \quad M=3\times7\times11\times13\times29\times31\times41\times43\times101\times151\times211\times281\times701
$$
contains no totients other than 10.
Theorem 2 (next section) implies that for almost all such $a$, the question
may be answered in the affirmative.
\bigskip
{\bf 3. A negative result}
\bigskip

\proclaim{Theorem 2} For any $\varepsilon>0$ there exist such $m$ that
at least $(1-\varepsilon)m$ residue classes $a\pmod{4m}$,
$0<a<4m$, $a\equiv 2\pmod4$ are totient-free.
\endproclaim

\proclaim{Corollary} The union of all totient-free residue classes
has density 3/4.
\endproclaim

\proclaim{Lemma 4} For any prime $r\ge5$ and for any $k=2,\dots,r-2$,
the number of distinct residues $x^k-x^{k-1}\pmod{r}$ with $(x,r)=1$
is less than $r-\sqrt{r/2}$.
\endproclaim
\noindent
Remark 1. The restriction $(x,r)=1$ is not essential as $0^k-0^{k-1}=
1^k-1^{k-1}$.

\noindent
Remark 2. Surely, the estimate of Lemma 4 is very weak, and it should be
$\le cr$, $c<1$. However, Lemma 4 is sufficient to prove Theorem 2.

\demo{Proof of Lemma 4} Let us consider the congruence
$$
x^k-x^{k-1}\equiv y^k-y^{k-1}\pmod{r},\quad 1\le x<r,\quad 1\le y<r,
\quad x\neq y.\tag3
$$
Let $y\equiv xz\pmod{r}$, $2\le z<r$. Any $z$ entails the unique solution
of (3) (namely, $x\equiv (z^{k-1}-1)/(z^k-1)$) if $z^{k-1}\not\equiv 1\pmod{r}$
and $z^k\not\equiv 1\pmod{r}$, otherwise $z$ does not entail any solutions.
So, the number of solutions of (3) is
$$
N=r-(r-1,k)-(r-1,k-1),
$$
since $(r-1,j)$ is the number of solutions to $z^j\equiv1\pmod{r}$.
Now $(r-1,k)$ and $(r-1,k-1)$ are coprime proper divisors of $r-1$.
Thus, their sum is at most $2+(r-1)/2$, so $N\ge(r-3)/2$.
If the number of distinct residues $x^k-x^{k-1}\pmod{r}$ with 
$(x,r)=1$ is $r-L$, then $L(L-1)\ge N$, hence $L^2\ge N+L>r/2$.
\qed\enddemo

Theorem 2 is equivalent to the following statement.

\proclaim{Theorem 2'} For any $\varepsilon>0$ there exist such odd 
$m$ that
for at least $(1-\varepsilon)m$ residues $a\pmod{m}$ the congruence (2)
does not have solutions with integers $k>0$ and $x$ with $(x,m)=1$.
\endproclaim

The equivalence of Theorems 2 and 2' follows directly from Lemma 3
and from the fact that the number of values of $a$ in (2) of the form
$p-1$ with $p$ a prime factor of $m$ is $O(\log m)$.

\proclaim{Lemma 5} For any $D\ge 1$ there are $\gg_D x/\log x$ primes $p\le x$
for which $D|(p-1)$ and no prime factor of $p-1$ exceeds $x^{9/20}$.
The result holds for $x$ sufficiently large depending on $D$.
\endproclaim

\demo{Proof} When $D=1$, this follows from the Theorem 1 of \cite{P}.
Since $D$ is fixed and $x\to\infty$, the general result follows by the same
method.
\qed\enddemo

\noindent
Remark 3.  The exponent $9/20$ in Lemma 5 is not the best possible
exponent.  For example, using the main theorem of \cite{F}, one
can replace $9/20$ with any number larger than $1/(2\sqrt{e})$.  However,
all we shall need below is an exponent smaller than $1/2$.
\demo{Proof of Theorem 2'} Let $p_1,\dots,p_I$ and $q_1,\dots,q_J$
be distinct odd primes such that
$$
\prod_i(1-1/p_i)<\varepsilon/4,\quad\prod_j(1-
1/q_j)<\varepsilon/4.\tag4
$$
Set $D=\lcm(p_1-1,\ldots, p_I-1, q_1-1, \ldots, q_J-1)$.
Let $y$ be a sufficiently large number and let
$r_1,\dots,r_L$ denote the primes $\le y$, different from all
$p_i$, $q_j$, for which each $r_l-1$ is divisible by $D$ and by
no prime $>y^{9/20}$.  By Lemma 5, $L\gg y/\log y$.
Take
$$
m=\prod_ip_i\prod_jq_j\prod_lr_l.
$$
By (4), the number of $a\pmod{m}$ satisfying
$$
\exists i\ a\equiv 1\pmod{p_i},\quad\exists j\ a\equiv -1\pmod{q_j}\tag5
$$
is at least $(1-\e/2)m$.
If $a$ satisfies (5) and $x$ is a solution of (2) with $(x,m)=1$
then $k\not\equiv 0\pmod{p_i-1}$ and $k\not\equiv 1\pmod{q_j-1}$,
therefore $k\not\equiv 0\pmod{r_l-1}$ and
$k\not\equiv 1\pmod{r_l-1}$ for all $l$.
For such $k$ we can estimate
the number of possible residues $a\pmod{r_l}$ by Lemma 4. Denote
$$n=\lcm(p_1-1,\dots,p_I-1,q_1-1,\dots,q_J-1,r_1-1,\dots,r_L-1)=
\lcm(r_1-1,\dots,r_L-1).$$
By construction,
$$
n \le \prod_{p\le y^{9/20}} p^{[\log y/\log p]} \le \exp\{y^{9/20}\log y\}. 
$$
By Lemma 4,
for any $k=1,\dots,n$ such that for each $l$,
$k\not\equiv 0\pmod{r_l-1}$ and
$k\not\equiv 1\pmod{r_l-1}$, the number of $a\pmod m$
for which there exists $x$ with $(x,m)=1$ satisfying (2)
does not exceed
$$
m\prod_l(1-1/\sqrt{2r_l})<m\exp(-L/\sqrt{2y}).
$$
Thus, the number of $a$ satisfying (5) for which a solution of (2) with
$(x,m)=1$ exists 
is less than $mn\exp(-L/\sqrt{2y}) \le \e m/2$ if $y$ is large enough.
\qed\enddemo

{\bf 4. A positive result}
\bigskip

\proclaim{Theorem 3} The set of all odd numbers $m$ such that for any $s\ge1$
and for any even $a$ the residue class $a\pmod{2^sm}$ contains
infinitely many totients, has a positive lower density.
\endproclaim

Call an odd number $m$ ``good'' if for any $a$ the congruence (2)
has a solution with positive integers $k$ and $(x,m)=1$. Theorem 3 has an
equivalent form:

\proclaim{Theorem 3'} The set of all good odd numbers has a positive 
lower density.
\endproclaim

\proclaim{Lemma 6 }  Suppose $f(x,y)$ is a polynomial absolutely
irreducible modulo $p$.  Then the number $N$ of solutions modulo $p$ of
$f(x,y) \equiv 0\pmod{p}$ satisfies
$$
|N-(p+1)| \le (d-1)(d-2)\sqrt{p} + d,
$$
where $d$ is the total degree of $f$.
\endproclaim
\demo{Proof}
In the case that $f$ is non-singular over $\overline{\Bbb F}_p$,
we use Weil's theorem.  The extra $d$ on the right of the inequality
is an upper estimate for the number of solutions ``at infinity".
If $f$ is singular, we use the principal result of Leep and Yeomans \cite{LY}.
\qed\enddemo

\proclaim{Lemma 7}  Suppose $p$ is a prime and $L,a,s,t$ are positive integers
with $(as,p)=1$.  Then the polynomial 
$$
f(x,y) = y^L(1-x^s)-ax^{t}
$$
is absolutely irreducible modulo $p$.
\endproclaim

\demo{Proof}
If $f(x,y)$ is reducible over $\overline{\Bbb F}_p$, then
$$
h(y) = y^L - \frac{ax^{t}}{1-x^s}
$$
is reducible over the field $k=\overline{\Bbb F}_p(x)$.  
By the criterion of Capelli and R\'edei (see Theorem 21 in \cite{S}),
this forces the existence of some $b$ in $k$ such that 
$ax^t/(1-x^s)=b^q$ for some prime $q$ dividing $L$, or
$ax^t/(1-x^s)=-4b^4$, in which case 4 divides $L$.  However,
since $s$ is coprime to $p$, $1-x$ divides $1-x^s$ to just the first
power, so neither possibility can occur.  
\qed\enddemo

\noindent
Remark 4.  It is also possible to give a direct proof of Lemma 7.
Over $\bar k$
we have the factorization
$$
h(y) = (y-r_1z) \cdots (y-r_L z),
$$
where each $r_i\in\overline{\Bbb F}_p$ satisfies $r_i^L=1$,
$z\in \bar k$, and $z^L=ax^t/(1-x^s)$.  Since $h$ is reducible over $k$, 
there exists a product
$$
(y-r_{i_1}z) \cdots (y-r_{i_j}z) \in k[y],
$$
where $j<L$.  In particular, the constant coefficient lies in $k$, whence
$z^j \in k$.  If $m$ is the smallest positive integer with $z^m\in k$,
then we have $m|L$, $m<L$.  Writing $u(x)=z^m$, we have
$$
u(x)^{L/m} = \frac{ax^{t}}{1-x^s}.
$$
As $1-x$ divides $1-x^s$ to just the first power, this equation 
is clearly impossible.

\proclaim{Lemma 8} There is a number $p_0$ such that for any prime $p>p_0$,
any positive integers $L\le p^{1/10}$ and $l\le L$ and any integer $a$ the
congruence (2) has a solution with $m=p$, $k\equiv l\pmod{L}$ and
$(x,p)=1$.
\endproclaim

\demo{Proof}  We may assume $a\not\equiv0\pmod{p}$.
To prove the lemma, it is enough to show the
existence of a solution $y$ of the congruence
$$
y^L(1-g)\equiv ag^l\pmod{p}\tag6
$$
with a primitive root $g$.
Indeed, we can let $x\equiv g^{-1}\pmod{p}$ and $k=l-uL$, where
$u$ is such that $y\equiv g^u\pmod{p}$.
We show a solution $y,g$ to (6) exists by estimating the number of solutions
of
$$
y^L(1-z^s)\equiv az^{sl}\pmod{p},\tag7
$$
where $s$ is a square-free divisor of $p-1$, and using inclusion-exclusion. 
By Lemma 7,
the polynomial $y^L(1-z^s)-az^{sl}$ is absolutely irreducible.
For a square-free
divisor $s$ of $p-1$, let $N_s$ be the number of solutions of (7).
For $s\le p^{1/5}$ we apply Lemma 6 and for larger $s$ we use the
trivial bound $N_s \le pL$.  Write $N_s=p+E_s$.  By inclusion-exclusion,
the number of solutions of (6) with a primitive root $g$ is
$$
\split
N &=  \sum_{s|p-1} \frac{\mu(s) N_s}s
\\
&\ge  p\prod_{\Sb q|p-1 \\ q\text{ prime} \endSb} (1-1/q) - \sum_{s|p-1}
\frac{|E_s|}{s} \\
&\ge \phi(p-1) - \sum_{\Sb s\le p^{1/5} \\ s|p-1 \endSb} (L+sl)^2\sqrt{p}/s
 - \sum_{\Sb s>p^{1/5} \\ s|p-1 \endSb} p^{9/10} \\
&\ge \tfrac12 \phi(p-1)
\endsplit
$$
provided $p$ is sufficiently large.
\qed\enddemo

\proclaim{Corollary}
Suppose $p_1<p_2<\cdots < p_r$ are odd primes larger than $p_0$, $m=p_1\cdots p_r$
and for any $j\ge 2$
$$
(p_j-1,\lcm(p_i-1:1\le i<j))\le p_j^{1/10}.
$$
Then $m$ is good.
\endproclaim

\demo{Proof}
Let $a$ be arbitrary.
Set $n_j=\lcm(p_i-1:1\le i< j)$ and $P_j=p_1\cdots p_j$
for each $j$.  We construct numbers 
$x_j,k_j$ inductively as follows.  Choose $x_1,k_1$ so that $(x_1,p_1)=1$ and
$x_1^{k_1}-x_1^{k_1-1} \equiv a \pmod{p_1}$.  For $j=2,\cdots,r$, Lemma
8 implies the existence of numbers $x_j,k_j$ for which $(x_j,P_j)=1$,
$x_j \equiv x_{j-1} \pmod{P_{j-1}}$, $k_j \equiv k_{j-1} \pmod{n_{j}}$ and
$x_j^{k_j}-x_j^{k_j-1}\equiv a\pmod{P_j}$.  The pair $(x_r,k_r)$ satisfies
(2) with $(x_r,m)=1$.
\qed\enddemo

Call an odd number $m$ ``forbidden" if $m=p_1\dots p_j$
where $p_1\le\dots\le p_j$ are primes and
$$
(p_j-1,\lcm(p_i-1:1\le i< j))>p_j^{1/10}.
$$

\proclaim{Lemma 9} The number of forbidden
numbers in $(x,2x]$ is $O(x/\log^5 x)$.
\endproclaim

Theorem 3' follows easily from Lemma 9.  Take some
$P\ge p_0$. Then for $x\ge 2P$ there are $\gg x/\log P$ positive
integers without prime factors $\le P$. If $m$ in $(x,2x]$ is not good,
the Corollary to Lemma implies $m$ is
divisible by a forbidden number
$>P^2$. By Lemma 9, there are $\ll x/\log^4 P$ such numbers. Therefore,
for sufficiently large $P$ and $x\ge 2P$ we get $\gg x/\log P$ good
numbers not exceeding $x$.

\demo{Proof of Lemma 9}  There is a constant $c>0$ so that whenever $n\ge 10$,
the number of divisors of $n$ is $\le n^{c/\log\log n}$.
By standard estimates from the distribution
of ``smooth'' numbers (see \cite{HT}), the number of
integers in $(x,2x]$ with all prime factors $\le x^{20c/\log\log x}$
is $O(x/\log^5 x)$.  Thus, we have to estimate the number $N$
of forbidden integers $m\in(x,2x]$ such that $p_j>x^{20c/\log\log x}$.
Denoting $l=m/p_j$, $n=\lcm(p_i-1:1\le i<j)=\lcm(p-1:p|l)$, we have
$$
(p_j-1,n)>x^{2c/\log\log x}.
$$
For fixed $l$ there are at most $x^{c/\log\log x}$ divisors of $n$, and for any
$d|n$ there are at most $2x/(dl)$ numbers $p_j>1$ for which $lp_j\le 2x$
and $p_j\equiv 1\pmod{d}$.  Summing over all divisors
$d>x^{2c/\log\log x}$,
we find that $l$ generates at most
$$
\sum_d 2x/(dl)<\sum_d 2x/(lx^{2c/\log\log x})\le 2x/(lx^{c/\log\log x})
$$
forbidden numbers. Further, taking the sum over $l$, we obtain the required
inequality $N \ll x/\log^5 x$.
\qed\enddemo

\enddocument